\documentclass[11pt,a4paper]{amsart}

\usepackage{amsfonts, amsmath, amssymb, amsthm}
\theoremstyle{plain}

\theoremstyle{definition}

\theoremstyle{remark}

\usepackage{txfonts}
\usepackage{mathbbol}
\usepackage{ae}

\usepackage{color}
\usepackage{graphicx}
\usepackage{url}
\usepackage[abs]{overpic}
\usepackage{subfigure}
\usepackage{booktabs}
\usepackage{psfrag}
\usepackage{paralist}
\usepackage[a4paper]{anysize}

\newcommand\RR{{\mathbb R}}

\renewcommand{\phi}{\varphi}

\graphicspath{{images/}}

\begin{document}

\title{Geometric Reasoning with \texttt{polymake}}
\author[Gawrilow \& Joswig]{Ewgenij Gawrilow and Michael Joswig}
\address{Ewgenij Gawrilow, Institut f\"ur Mathematik, MA 6-1, TU Berlin, 10623 Berlin, Germany}
\email{gawrilow@math.tu-berlin.de} 
\address{Michael Joswig, Fachbereich Mathematik, AG~7, TU Darmstadt, 64289 Darmstadt, Germany}
\email{joswig@mathematik.tu-darmstadt.de}
\thanks{The second author is partially supported by Deutsche Forschungsgemeinschaft, DFG Research Group
  ``Polyhedral Surfaces.''}
\date{\today}

\begin{abstract}
  The mathematical software system \texttt{polymake} provides a wide range of functions for convex
  polytopes, simplicial complexes, and other objects.  A large part of this paper is dedicated to a
  tutorial which exemplifies the usage.  Later sections include a survey of research results
  obtained with the help of \texttt{polymake} so far and a short description of the technical
  background.
\end{abstract}

\maketitle

\section{Introduction}
\noindent
The computer has been described as \emph{the} mathematical machine.  Therefore, it is nothing but
natural to let the computer meet challenges in the area where its roots are: mathematics.  In fact,
the last two decades of the 20th century saw the ever faster evolution of many mathematical software
systems, general and specialized.  Clearly, there are areas of mathematics which are more apt to
such an approach than others, but today there is none which the computer could not contribute to.

In particular, polytope theory which lies in between applied fields, such as optimization, and more
pure mathematics, including commutative algebra and toric algebraic geometry, invites to write
software.  When the \texttt{polymake} project started in 1996, there were already a number of
systems around which could deal with polytopes in one way or another, e.g., convex hull codes such
as \texttt{cdd}~\cite{cdd}, \texttt{lrs}~\cite{lrs}, \texttt{porta}~\cite{porta}, and
\texttt{qhull}~\cite{qhull}, but also visualization classics like \texttt{Geomview}~\cite{geomview}.
The basic idea to \texttt{polymake} was ---and still is--- to make interfaces between any of these
programs and to continue building further on top of the combined functionality.  At the same time
the gory technical details which help to accomplish such a thing should be entirely hidden from the
user who does not want to know about it.  On the outside \texttt{polymake} behaves somewhat similar
to an expert system for polytopes: The user once describes a polytope in one of several natural ways
and afterwards he or she can issue requests to the system to compute derived properties.  In
particular, there is no need to program in order to work with the system.  On the other hand, for
those who do want to program in order to extend the functionality even further, \texttt{polymake}
offers a variety of ways to do so.

The modular design later, since version 2.0, allowed \texttt{polymake} to treat other mathematical
objects in the same way.  Mostly guided by the research of the second author the system was
augmented by the \texttt{TOPAZ} application which deals with finite simplicial complexes.  Because
of the connections between polytope theory and combinatorial topology both parts of the system now
benefit from each other.

We explain the organization of the text.  It begins with a quite long tutorial which should give an
idea of how the system can be used.  The examples are deliberately chosen to be small enough that it
is possible to verify all claims while reading.  What follows is an overall description of the key
algorithms and methods which are available in the current version 2.1.0.  The subsequent
Section~\ref{sec:survey} contains several brief paragraphs dedicated to research in mathematics that
was facilitated by \texttt{polymake}.

Previous reports on the \texttt{polymake} system include the two
papers~\cite{MR1785292,SoCG01:polymake}; by now they are partially outdated.  \texttt{polymake} is
open source software which can be downloaded from \url{http://www.math.tu-berlin.de/polymake} for
free.

\subsection*{Acknowledgments}

Over the time many people made small and large contributions to the code.  Most notably, Thilo
Schr\"oder and Nikolaus Witte are members of the development team since 2002.

Partial funding for the initial period of 1996--1997 of the \texttt{polymake} project came from the
German-Israeli Foundation for Scientific Research and Development, grant I-0309-146.06/93 of
G\"unter M. Ziegler.  Later the Deutsche Forschungsgemeinschaft (DFG) partially supported
\texttt{polymake} within projects of the Sonderforschungsbereich 288 ``Differentialgeometrie und
Quantenphysik'' and the DFG-For\-schungs\-zen\-trum Matheon.

\section{A Tutorial}
\noindent
This tutorial tries to give a first idea about \texttt{polymake}'s features by taking a look at a
few small examples.  We focus on computations with convex polytopes.  For definitions, more
explanations, and pointers to the literature see the subsequent Section~\ref{subsec:polytopes}.

The text contains commands to be given to the \texttt{polymake} system (preceded by a '$>$') along
with the output.  As the environment a standard UNIX shell, such as \texttt{bash}, is assumed.
Commands and output are displayed in \texttt{typewriter type}.

Most of the images shown are produced via \texttt{polymake}'s interface to
\texttt{JavaView}~\cite{javaview} which is fully interactive.

\subsection{A very simple example: the $3$-cube}

Suppose you have a finite set of points in the Euclidean space $R^d$.  Their convex hull is a
polytope $P$.  Now you want to know how many facets $P$ has.  As an example let the set $S$ consist
of the points $(0,0,0)$, $(0,0,1)$, $(0,1,0)$, $(0,1,1)$, $(1,0,0)$, $(1,0,1)$, $(1,1,0)$, $(1,1,1)$
in $R^3$.  Clearly, $S$ is the set of vertices of a cube.  Employing your favorite word processor
produce an ASCII text file, say \texttt{cube.poly}, containing precisely the following information.

\begin{verbatim}
POINTS
1 0 0 0
1 0 0 1
1 0 1 0
1 0 1 1
1 1 0 0
1 1 0 1
1 1 1 0
1 1 1 1
\end{verbatim}

We have the keyword \texttt{POINTS} followed by eight lines, where each line contains the
coordinates of one point.  The coordinates in $\RR^3$ are preceded by an additional $1$ in the first
column; this is due to the fact that \texttt{polymake} works with \emph{homogeneous coordinates}.
The solution of the initial problem can now be performed by \texttt{polymake} as follows.
Additionally, we want to know whether $P$ is simple, that is, whether each vertex is contained in
exactly $3$ facets (since $\dim P=3$).

\begin{verbatim}
> polymake cube.poly N_FACETS SIMPLE
\end{verbatim}

While \texttt{polymake} is searching for the answer, let us recall what \texttt{polymake} has to do
in order to obtain the solution.  It has to compute the convex hull of the given point set in terms
of an explicit list of the facet describing inequalities.  Then they have to be counted, which is,
of course, the easy part.  Nonetheless, \texttt{polymake} decides on its own what has to be done,
before the final ---admittedly trivial--- task can be performed.  Checking simplicity requires to
look at the vertex facet incidences.  In the meantime, \texttt{polymake} is done.  Here is the
answer.

\begin{verbatim}
N_FACETS
6

SIMPLE
1
\end{verbatim}

Simplicity is a boolean property.  So the answer is either ``yes'' or ``no'', encoded as $1$ or $0$,
respectively.  The output says that $P$ is, indeed, simple.

Depending on the individual configuration \texttt{polymake} chooses one of several convex hull
computing algorithms.  In the previous example \texttt{polymake} might have used the double
description method from Fukuda's cdd package~\cite{cdd}.  It is possible to explicitly specify other
methods.

As a matter of fact, \texttt{polymake} knows quite a bit about standard constructions of polytopes.
So you do have to type in your $3$-cube example.  You can use the following command instead.  The
trailing argument $0$ indicates a cube with $0/1$-coordinates.

\begin{verbatim}
> cube cube.poly 3 0
\end{verbatim}

\subsection{Visualizing a Random Polytope}

But let us now try something else.  How does a typical polytope look like?  To be more precise: Show
me an instance of the convex hull of $20$ randomly distributed points on the unit sphere in $R^3$.
This requires one command to produce a \texttt{polymake} description of such a polytope and a second
one to trigger the visualization.  Again there is an implicit convex hull computation going on
behind the scenes.  On the way a complete combinatorial description of the polytope is obtained.

\begin{verbatim}
> rand_sphere random.poly 3 20
> polymake random.poly VISUAL
\end{verbatim}

\texttt{polymake}'s standard tool for interactive visualization is \texttt{JavaView} by Polthier and
others~\cite{javaview} For instance, it allows you to rotate or zoom into your polytope.  Here is a
sequence of snapshots.

\begin{figure}[htbp]\centering
  \includegraphics[height=.23\textwidth]{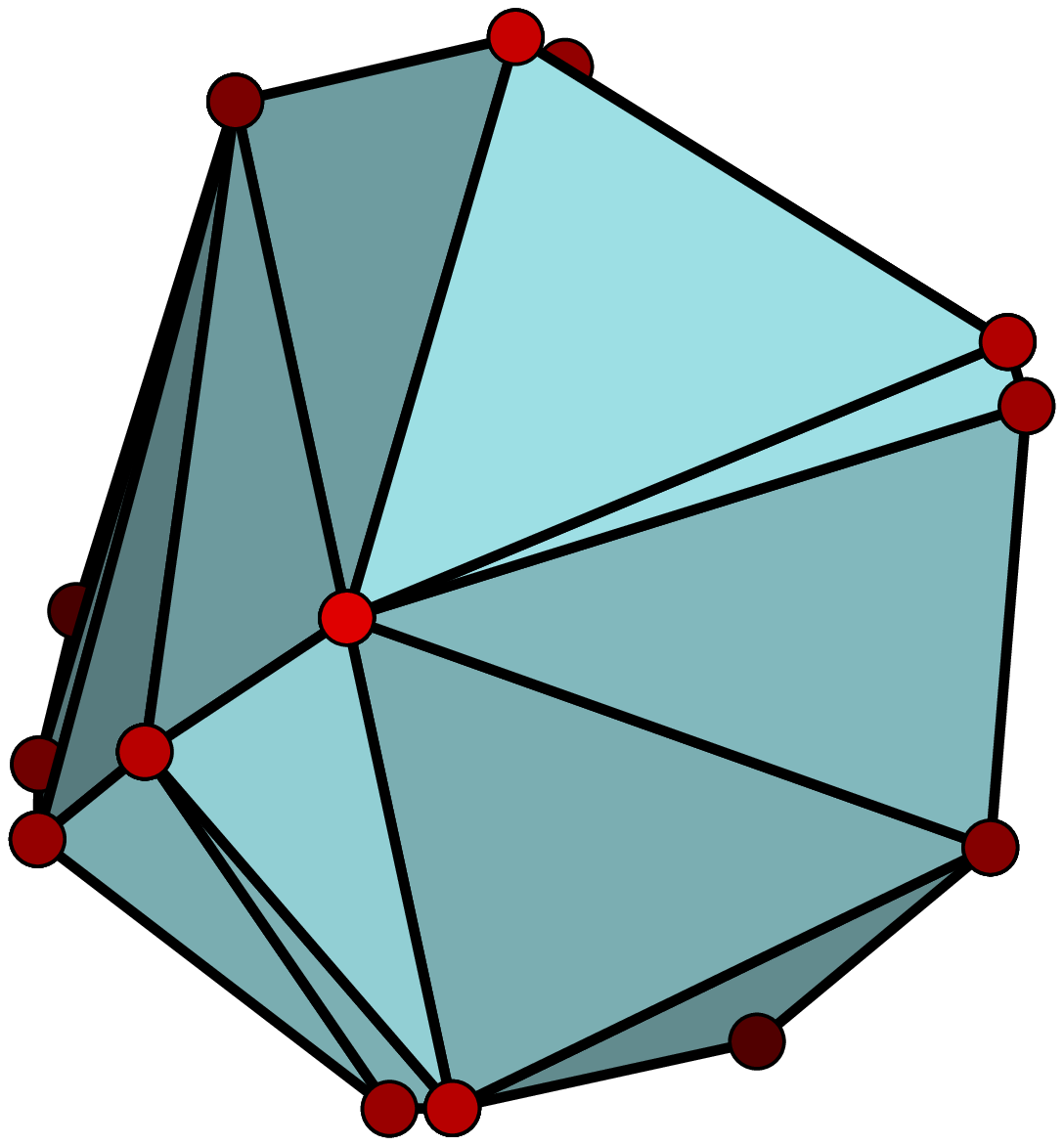}\quad
  \includegraphics[height=.23\textwidth]{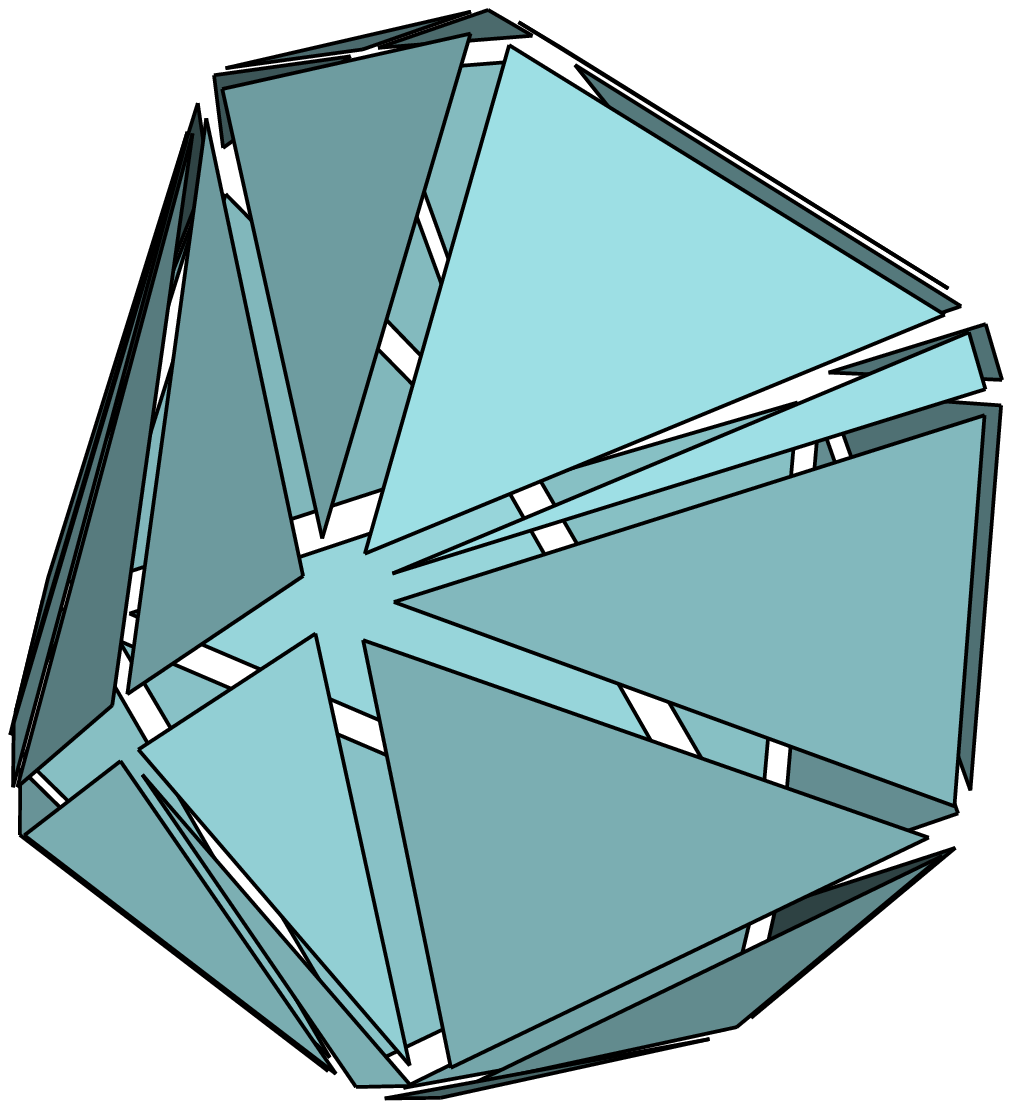}\quad
  \includegraphics[angle=90,height=.23\textwidth]{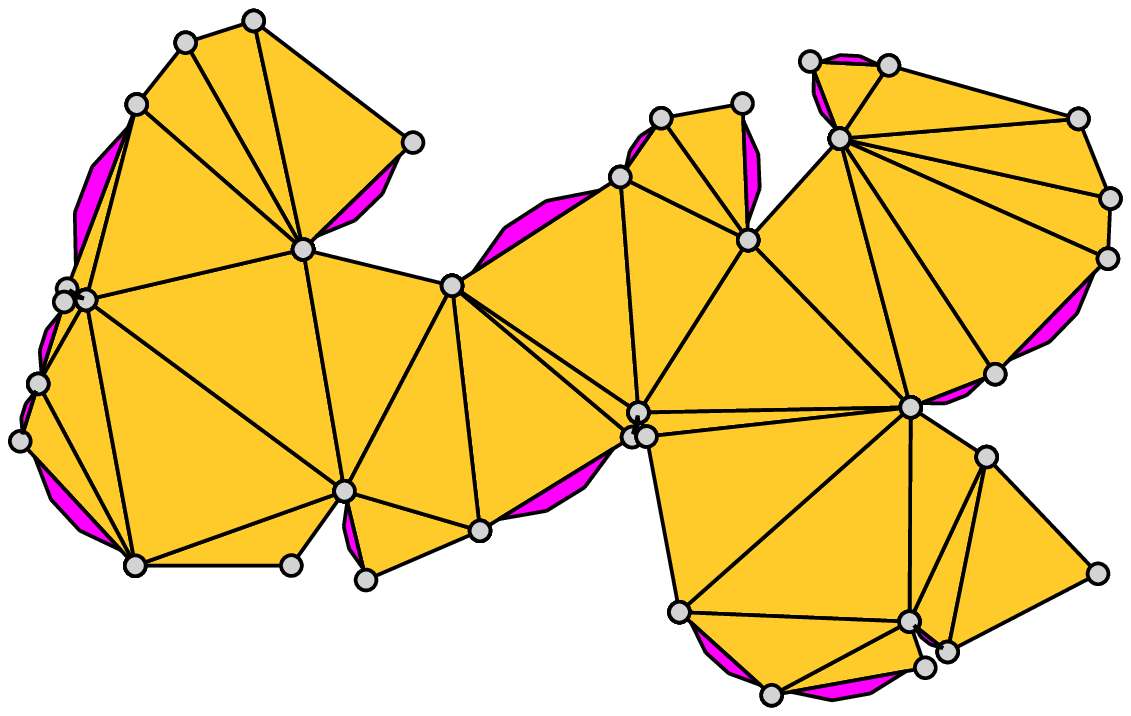}
  \caption{Interactive visualization of a random polytope with \texttt{JavaView}: Three
    snapshots.\label{fig:interactive}}
\end{figure}

\subsection{Linear programming}

Polytopes most naturally appear as sets of feasible solutions of linear programs.  Consider the
following example.

\begin{figure}[htbp]\centering
  \begin{minipage}{.5\textwidth}
    \begin{align*}
      \text{Max}&\text{imize}\\
      & x_1 + x_2 + x_3\\
      \text{sub}&\text{ject to}\\
      & 0\le x_1, x_2, x_3 \le 1\\
      & x_1 + x_2 + x_3 \le 5/2\\
      & x_1 - 17x_2 \le 8
    \end{align*}
  \end{minipage}
  \hfill
  \begin{minipage}{.3\textwidth}
    \includegraphics[width=\textwidth]{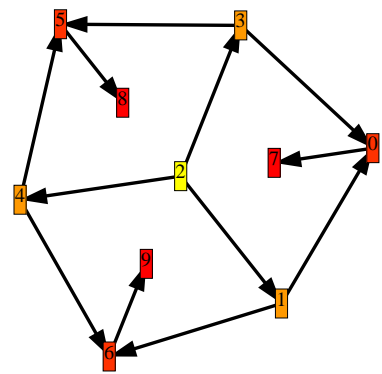}
  \end{minipage}
  \caption{Small linear program and a visualization.\label{fig:small-linear}}
\end{figure}

A linear inequality $a_0 + a_1 x_1 + \dots a_d x_d \ge0$ is encoded as the inward pointing normal
vector to the corresponding affine hyperplane (of suitable length).  This is what \texttt{polymake}
uses: the former inequality is represented as the vector $(a_0,a_1,\dots,a_d)$.  Our linear program
in \texttt{polymake}'s format looks as given below.  Note that we do not decide beforehand, whether
we want to minimize or maximize.

\begin{verbatim}
LINEAR_OBJECTIVE
 0    1  1  1

INEQUALITIES
 0    1  0  0
 0    0  1  0
 0    0  0  1
 1   -1  0  0
 1    0 -1  0
 1    0  0 -1
5/2  -1 -1 -1
 8   -1 17  0
\end{verbatim}

People working in optimization usually prefer other file formats (such as \texttt{CPLEX}'s LP file
format), where it is also possible to keep the names of the variables.  \texttt{polymake} is aware
of this.  LP format files can be converted by calling the command \texttt{lp2poly}.
      
It is not difficult to determine the polytope forming the solution space.  Assume the file
\texttt{linear\_{}program.poly} contains the description given above.

\begin{verbatim}
> polymake linear_program.poly MAXIMAL_VALUE
MAXIMAL_VALUE
5/2
\end{verbatim}

This is the kind of behavior one would expect from a linear solver.  Of course, usually it is also interesting to
obtain a point attaining the maximum.  And, in fact, \texttt{polymake} calls \texttt{cdd}'s
implementation of the Simplex Method (with exact rational arithmetic).
With \texttt{polymake} you can go one step further.  It visualizes for you the polytope with
directed edges.  Instead of relying on the interactive \texttt{JavaView} visualization this time we
produce postscript output directly; see Figure~\ref{fig:small-linear}.

\begin{verbatim}
> polymake linear_program.poly postscript \
           "VISUAL_GRAPH->DIRECTED_GRAPH->VERTEX_COLORS"
\end{verbatim}

The directed edges (whose orientation is induced by the given linear objective function) are drawn as
arrows.  The colors of the vertices indicate the level of height with respect to the objective
function.  The numbering of the vertices corresponds to the order of the vertices in the file.

It is an essential feature of polytope theory that it is possible to define a polytope in one of two
equivalent ways: either as the convex hull points or as the intersection of half-spaces.  This is
reflected in \texttt{polymake}: All functions can directly be applied to all polytopes regardless of
how they were initially described.

\subsection{Polytopes From the Combinatorial Point of View}

Suppose now you are not interested in a particular coordinate representation of a polytope.  But instead you want
to focus on the combinatorial properties only.  \texttt{polymake} supports this point of view, too.  You can specify a
polytope in terms of its vertex-facet-incidence matrix.  For each facet you have a line with a list of the vertices
contained in that facet.  The vertices are specified by numbers.  They are numbered consecutively starting from $0$.
In each row the vertices are listed in ascending order.  The following is a valid \texttt{polymake} description of a
square.

\begin{verbatim}
VERTICES_IN_FACETS
{0 1}
{1 2}
{2 3}
{0 3}
\end{verbatim}

Note that in this situation \texttt{polymake} assumes that you actually specified a polytope in this
way.  Verifying realizability, albeit algorithmically possible (e.g., by means of quantifier
elimination), is beyond \texttt{polymake}'s capabilities.

The dimension of a polytope, i.e., the dimension of its affine hull, is an intrinsic property.  It
does not depend on the coordinate representation.  The vertex-facet-incidence matrix suffices for
\texttt{polymake} to compute the dimension.  Assume the data above was stored in a file named
\texttt{square.poly}.
      
\begin{verbatim}
> polymake square.poly DIM
DIM
2
\end{verbatim}

Often one wants to construct new polytopes from old ones.  So suppose we need a prism over a
triangle.  This can be constructed as the \emph{wedge} of a square over an arbitrary facet (e.g.,
the first facet, which is numbered 0).  This is a simple polytope, but it is not simplicial.  The
program option \texttt{-noc} means ``no coordinates'', and this limits the produced output to a
purely combinatorial description.

\begin{verbatim}
> wedge prism.poly square.poly 0 -noc
> polymake prism.poly SIMPLE SIMPLICIAL
SIMPLE
1

SIMPLICIAL
0
\end{verbatim}

\subsection{Checking for Combinatorial Equivalence}

Suppose the file \texttt{Sharir.poly} contains the following inequality description of a $3$-dimensional polytope:

\begin{verbatim}
INEQUALITIES
25 -2 -25 10
-2 25 2 10
25 -2 25 10
-2 25 -2 10
0 0 -1 -1
2 0 -1 1
\end{verbatim}

It turns out that the polytope defined this way is, in fact, combinatorially equivalent to the $3$-dimensional
cube that we studied above.  \texttt{polymake} can give a proof of this as shown below; internally
McKay's \texttt{nauty} is called~\cite{nauty}.
\begin{verbatim}
> polymake -v check_iso Sharir.poly cube.poly
[fixing partition]
(2 6)(3 7)(8 11)(9 10)
level 3:  10 orbits; 3 fixed; index 2
(1 3)(5 6)(8 13)(9 12)
level 2:  6 orbits; 1 fixed; index 3
(0 1)(2 3)(4 5)(6 7)(12 13)
level 1:  2 orbits; 0 fixed; index 8
2 orbits; grpsize=48; 3 gens; 10 nodes; maxlev=4
tctotal=20; canupdates=1; cpu time = 0.00 seconds
[fixing partition]
(2 4)(3 5)(10 12)(11 13)
level 3:  10 orbits; 2 fixed; index 2
(1 2)(5 6)(8 10)(9 11)
level 2:  6 orbits; 1 fixed; index 3
(0 1)(2 3)(4 5)(6 7)(8 9)
level 1:  2 orbits; 0 fixed; index 8
2 orbits; grpsize=48; 3 gens; 10 nodes; maxlev=4
tctotal=20; canupdates=1; cpu time = 0.00 seconds
h and h' are identical.
 0-0 1-1 2-3 3-2 4-7 5-6 6-5 7-4 8-10 9-11 10-13 11-12 12-9 13-8
check_iso
1
\end{verbatim}

What you see is output for a nauty computation which gives you some information about the orbit
structure of the automorphism groups of the polytopes involved.  Essential is the fourth to last
line telling us that both polytopes are, in fact, combinatorially equivalent.  The third to last
line describes one particular isomorphism: For each vertex (one of the nodes numbered $0$ through
$7$) and each facet (one of the nodes numbered $8$ through $13$) of \texttt{Sharir.poly} the
corresponding vertex or facet of \texttt{cube.poly} is listed.  The final two lines are
\texttt{polymake}'s output indicating that the two polytopes are indeed combinatorially equivalent.
Omitting the \texttt{-v} flag in the command line would have resulted in these two output lines only.

\subsection{A More Detailed Look}

Let us compute the volume of a polytope.  And we want to know what \texttt{polymake} actually does
in order to obtain the result.  Adding the \texttt{-v} or \texttt{-vv} flag to the \texttt{polymake}
command line call asks for (very) verbose information.  We omit a couple of lines at the beginning
of the output, where the system tells about its rule base.

The output below corresponds to computing the \texttt{VOLUME} directly from the \texttt{INEQUALITIES}
description.  It looks different if you solved that linear optimization problem before.

\begin{verbatim}
> polymake -vv linear_program.poly VOLUME
polymake: reading rules from ...
polymake: minimum weight rule chain constructed in 0.054 sec.
polymake: applying rule cdd.convex_hull.dual: VERTICES, POINTED, FEASIBLE :
                                              FACETS | INEQUALITIES
polymake: applying rule BOUNDED : VERTICES | POINTS
polymake: applying rule PRECONDITION: BOUNDED ( default.volume: VOLUME :
                                                VERTICES, TRIANGULATION )
polymake: applying rule beneath_beyond.convex_hull.primal, default.triangulation:
 FACETS, AFFINE_HULL, VERTICES_IN_FACETS, DUAL_GRAPH, TRIANGULATION,
 ESSENTIALLY_GENERIC : VERTICES
polymake: applying rule default.volume: VOLUME : VERTICES, TRIANGULATION
VOLUME
47/48
\end{verbatim}

The program which finally produces the volume is very simple-minded: It takes any triangulation and
adds up the volumes of the simplices.  But before that, \texttt{polymake} does the following: From
the rule base, the system infers that it should first call \texttt{cdd} to obtain the vertices from
the inequalities via a (dual) convex hull computation.  Then it checks whether the input polyhedron
is bounded, that is, to check whether the volume is finite; otherwise \texttt{polymake} would abort
the computation.  As a third step \texttt{polymake} constructs a triangulation, which, in fact, is
obtained from calling a second convex hull code, which differs from \texttt{cdd}'s double
description method in that it additionally produces a triangulation.

An important feature of \texttt{polymake} is that all intermediate data which are computed are
stored into the polytope file.  Asking the program a second time for the same thing, or, for
something else which had been computed before, gives an immediate answer.  The result is read from
the file.

\begin{verbatim}
> polymake -vv linear_program.poly VOLUME
polymake: reading rules from ...
VOLUME
47/48
\end{verbatim}

\texttt{polymake} employs a special transaction model in order to assert the consistency of the
polytope data.  It relies on its rule base in the sense that rules whose execution terminates
properly are expected to produce valid output.  Moreover, \texttt{polymake} has a powerful error
recovery mechanism which automatically looks for an alternative way to find the answer to a user's
question in case of the failure of some external program.

\section{Applications Overview}
\noindent
There are several different kinds of mathematical objects which \texttt{polymake} can deal with,
most notably convex polytopes and finite simplicial complex.

As it was shown in the tutorial section above the system's behavior is driven by rules.  In fact,
the part of the system which takes care of applying rules to answer user requests is entirely
independent of the mathematical objects and the algorithms.  Each class of objects comes with its
own set of rules.  We survey the more technical aspects in Section~\ref{sec:design}.  Here we focus
on the mathematics.

\subsection{Convex Polytopes}
\label{subsec:polytopes}

A convex polytope is the convex hull of finitely many points in $\RR^d$; this is its
\emph{V-description}.  A basic result in this area says that this notion coincides with those
intersections of finitely many affine halfspaces in $\RR^d$ which are bounded
(\emph{H-description}).  This is sometimes referred to as the \emph{Main Theorem on convex
  polytopes}.  For an introduction to the theory, see Gr\"unbaum~\cite{Gruenbaum1967} or
Ziegler~\cite{Ziegler1995}.

Convex polytopes arise in many mathematical fields as diverse as linear and combinatorial
optimization, combinatorial topology, commutative algebra and algebraic geometry.  All these areas,
in their modern form, pursue algorithmic questions, and this is where \texttt{polymake} proved to be
useful.  The Section~\ref{sec:survey} discusses some of these applications in more detail.
  
In order to deal with polytopes algorithmically often a first step is to apply an effective version
of the ``Main Theorem''.  While a polytope may naturally be given in its H-description (such as in
linear programs) it is essential to also obtain a V-representation if one is interested in
combinatorial properties.  Algorithms which solve this problem are \emph{convex hull algorithms}.
Many such algorithms are known an implemented.  The running-time that a particular
algorithm/implementation requires is known to vary strongly with the class of polytopes which it is
applied to; see Avis, Bremner, and Seidel~\cite{MR1447243} and also~\cite{MR2011751}.  Therefore,
\texttt{polymake} offers three different convex hull algorithms for the user to choose.  There is
one which is built into the system and two more, \texttt{cdd}~\cite{cdd} and \texttt{lrs}~\cite{lrs},
respectively, which are accessible via interfaces.  Actually, there also interfaces to
\texttt{porta}~\cite{porta} and \texttt{qhull}~\cite{qhull} available, but these are disabled by
default.  For each call it is possible to specify which algorithm to choose; additionally, the
system can freely be configured to generally prefer one algorithm over another.

Convex polytopes have a metric geometry look as well as a combinatorial one. Both views are
supported by \texttt{polymake}.  What follows first is a list of metric properties which can be
computed with the software.
\begin{itemize}
\item Gale transformations
\item Steiner points
\item projective linear transformations
\item triangulations
\item Voronoi diagrams and Delaunay cell decompositions in arbitrary dimension
\end{itemize}
Combinatorial properties include:
\begin{itemize}
\item fast face lattice construction algorithm; due to Kaibel and Pfetsch~\cite{MR1927137}
\item $f$-vector, $h$-vector, flag-$f$-vector, and $cd$-index
\item various graph-theoretic properties of the vertex-edge graph
\item Altshuler determinant
\end{itemize}

\begin{figure}[htbp]\centering
  \includegraphics[width=.45\textwidth]{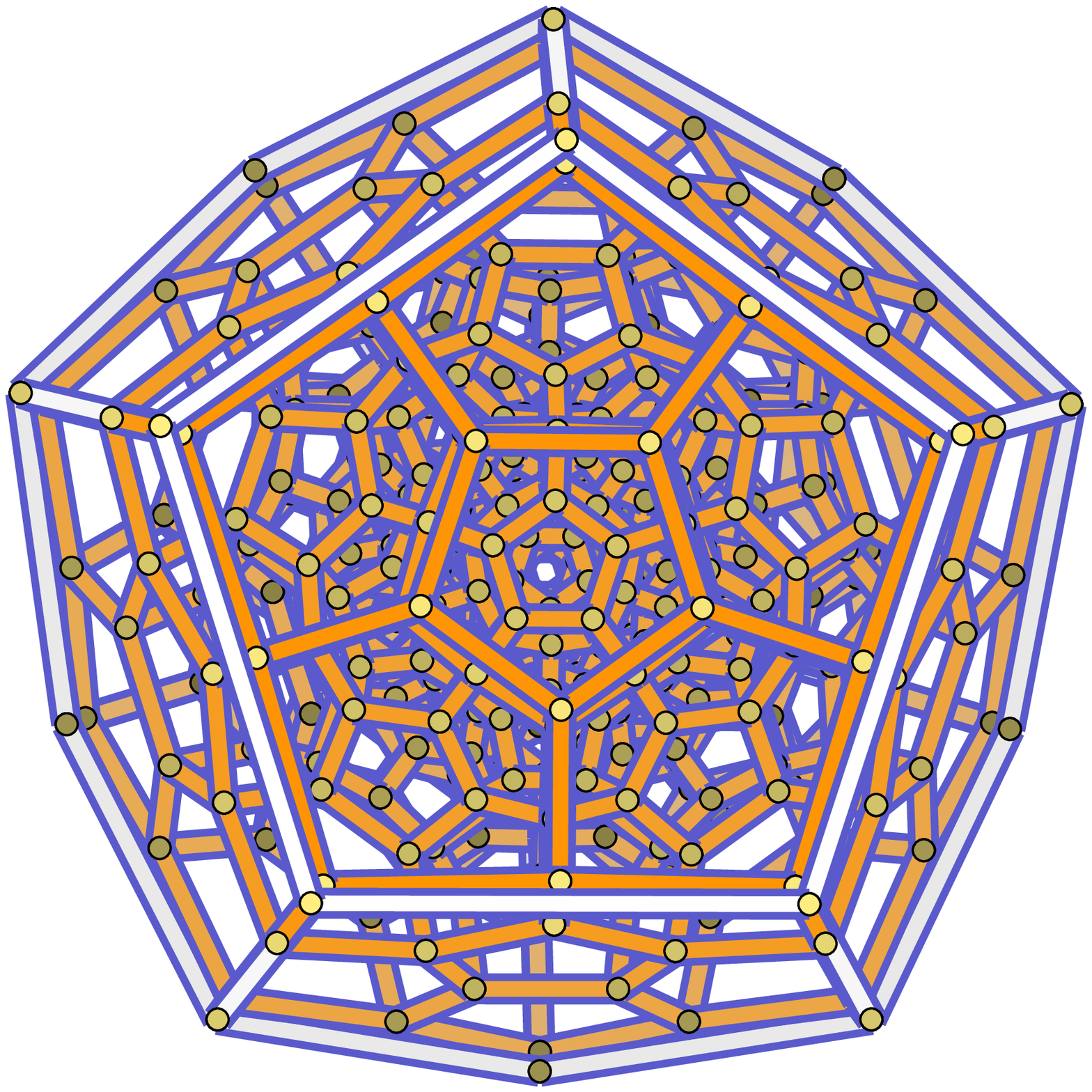}\qquad
  \includegraphics[width=.45\textwidth]{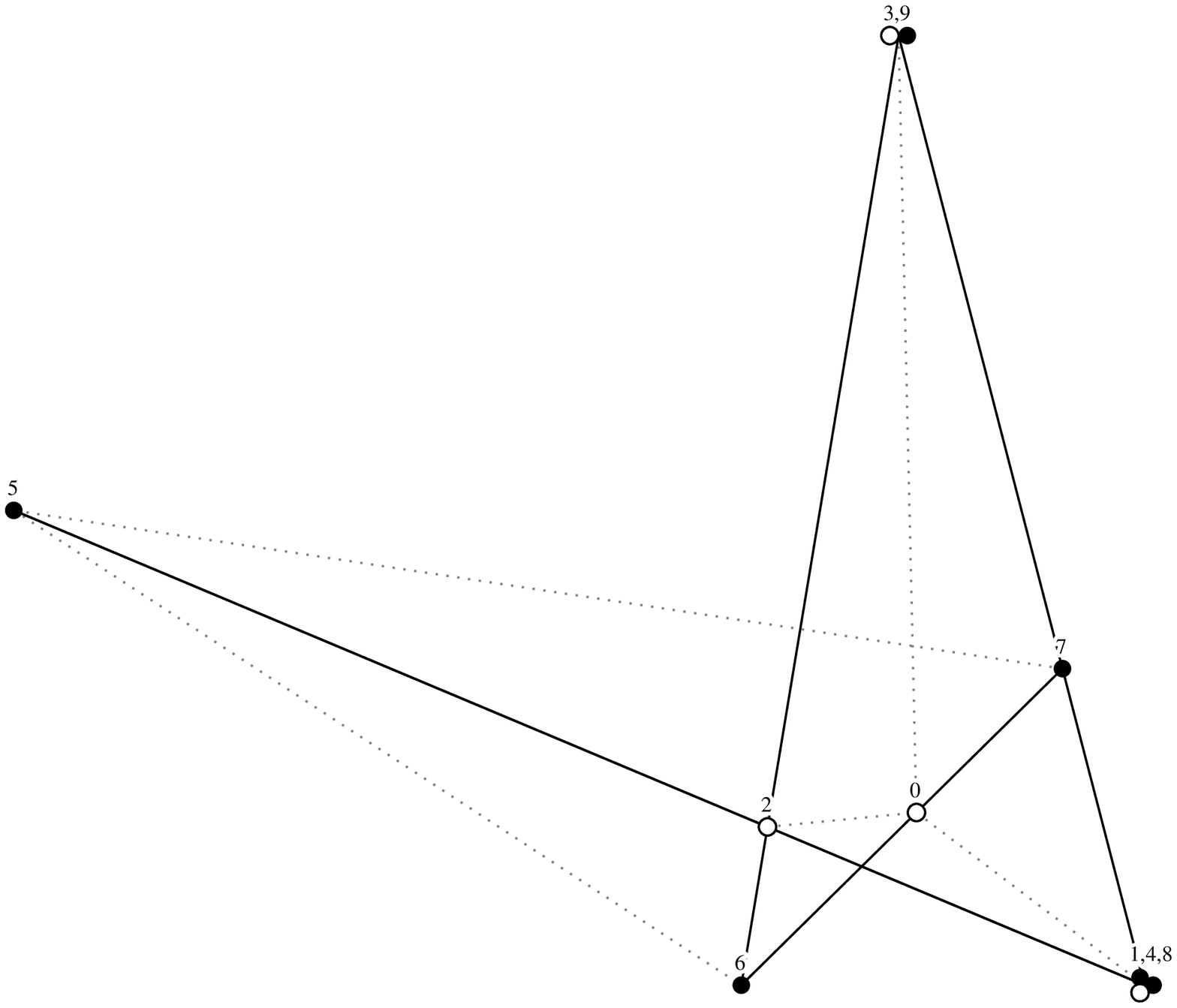}\qquad
  \caption{Schlegel diagram of the regular $120$-cell (left) and a Gale diagram of a random
    $6$-dimensional $01$-polytope with $10$ vertices (right).\label{fig:visualization}}
\end{figure}

In addition to these features there is a wide range of standard constructions and visualization
functions; e.g., see Figure~\ref{fig:visualization}.  There is also an interface to
\texttt{Geomview}~\cite{geomview}.

\subsection{Finite Simplicial Complexes}

Given a finite \emph{vertex set} $V$, a \emph{simplicial complex} on~$V$ is a subset of $2^V$ which
is closed with respect to taking subsets.  Simplicial complexes form a basic combinatorial concept
to capture properties of well-behaved topological spaces.  In particular, this way certain parts of
topology get within reach of effective methods.

A fundamental problem in topology is to decide whether two given spaces are \emph{homeomorphic},
that is, indistinguishable from the topological point of view, or not.  While it can be shown that
this is algorithmically impossible ---even for finite simplicial complexes representing
$4$-dimensional manifolds--- it remains a key task to compute algebraic (homotopy) invariants.

\texttt{polymake} offers the following:
\begin{itemize}
\item simplicial homology and cohomology with integer coefficients
\item cup and cap products
\item Stiefel-Whitney characteristic classes
\item intersection forms of $4$-manifolds
\item flip-heuristic by Bj\"orner and Lutz~\cite{MR2001h:57026} for detecting spheres
\end{itemize}

In particular, in view of a celebrated result of Freedman~\cite{MR84b:57006}, \texttt{polymake} is
able to solve the homeomorphism problem for combinatorial $4$-manifolds which are simply connected.
See the survey~\cite{math.AT/0401176} for some example computations.

\subsection{Extensions and Related Concepts}

The whole \texttt{polymake} system is extensible in several ways.  Besides adding new functionality
to the applications dealing with polytopes and simplicial complexes, it is possible to define
entirely new classes of objects with an entirely new set of rules.  For the more technical aspects
such an extension the reader is referred to Section~\ref{sec:design}.

Here we list features which are already built into the system but which go beyond standard
computations with polytopes or simplicial complexes.

\subsubsection{Tight Spans of Finite Metric Spaces}\label{subsubsec:tightspan}

\begin{figure}[htbp]
  \begin{minipage}{.57\textwidth}\small
    \renewcommand{\arraystretch}{0.9}
    \begin{tabular*}{\linewidth}{@{\extracolsep{\fill}}lrrrrrr@{}}\toprule
             & Athens & Berlin & Lisboa & London & Paris & Rome\\
      \midrule
      Athens &      0 &   1535 &   2815 &   1988 &  1827 & 1522\\
      Berlin &        &      0 &   1864 &    727 &   677 &  926\\
      Lisboa &        &        &      0 &   1442 &  1114 & 1672\\
      London &        &        &        &      0 &   280 & 1125\\
      Paris  &        &        &        &        &     0 &  870\\
      Rome   &        &        &        &        &       &    0\\
      \bottomrule
    \end{tabular*}
  \end{minipage}
  \hfill
  \begin{minipage}{.4\textwidth}
    \includegraphics[width=\textwidth]{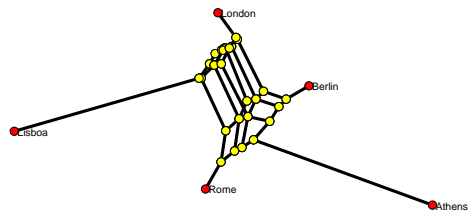}
  \end{minipage}
  \caption{Distances (in kilometers) among six European cities and a $3$-dimensional visualization
    of their tight span, which lives in~$\RR^6$.\label{fig:tightspan}}
\end{figure}

Every tree $T$ with non-negative weights on the edges defines a metric on the nodes of $T$.
Conversely, it is easy to reconstruct the tree from such a \emph{tree-like} metric.  The
\emph{phylogenetic problem} in computational biology boils down to the task to derive \emph{a}
sufficiently close tree from any given finite metric space.  It is obvious that sometimes there is
no tree at all which fits a given metric.  Dress et al.~\cite{MR97i:57002,MR2003g:54077} devised
\emph{tight spans} as geometric objects which can be assigned to any finite metric space and which
capture the deviation from a tree-like metric.  Since tight spans can be described as bounded
subcomplexes of unbounded polyhedra, \texttt{polymake}'s features can be exploited.  See
Figure~\ref{fig:tightspan} for an example.

Sturmfels and Yu~\cite{MR2097310} recently used \texttt{TOPCOM}~\cite{topcom} and \texttt{polymake}
to classify tight spans of metric spaces with at most six points.

\subsubsection{Curve Reconstruction}

\begin{figure}[htbp]\centering
  \includegraphics[width=.3\textwidth,clip=true]{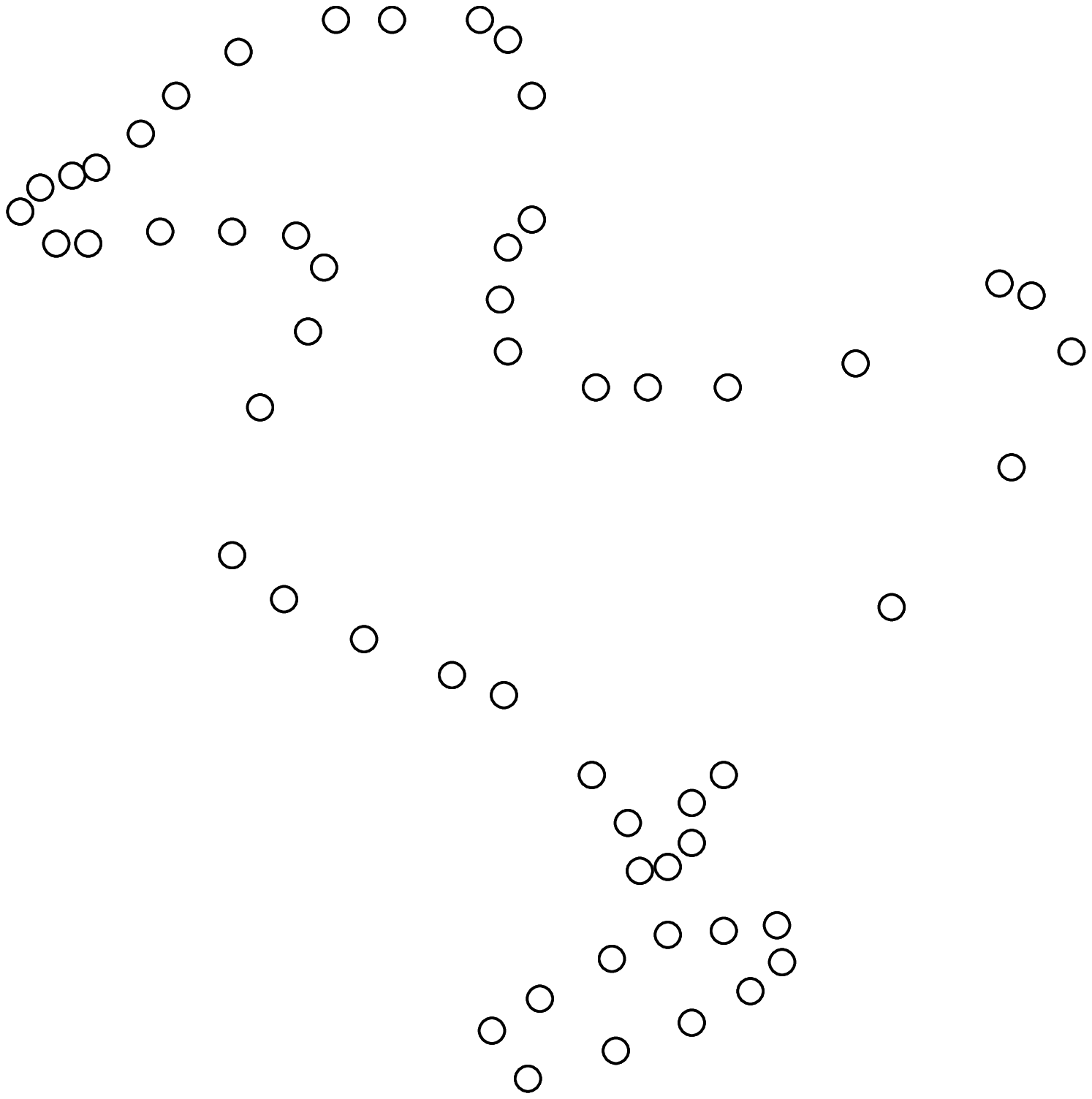}\hfill
  \includegraphics[width=.3\textwidth,clip=true]{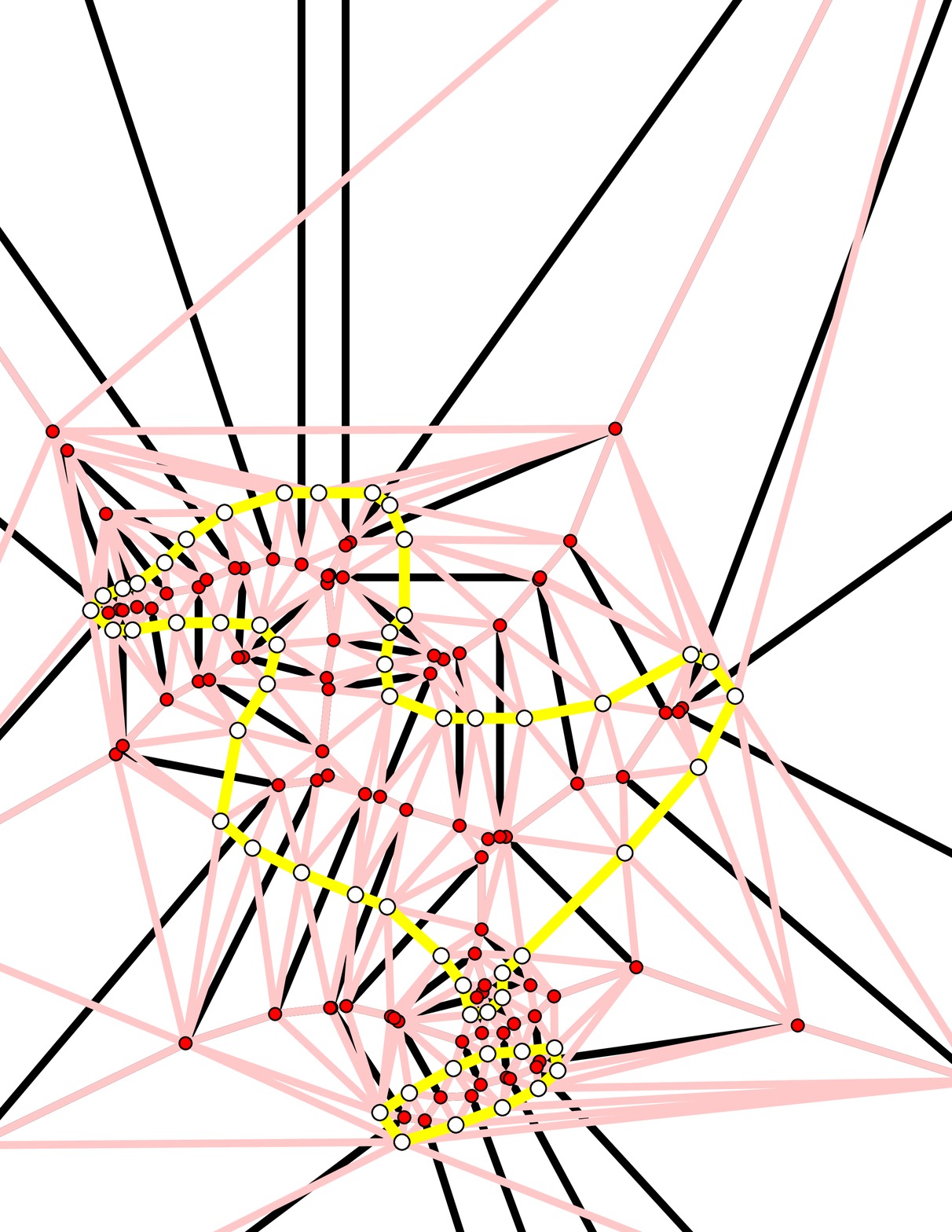}\hfill
  \includegraphics[width=.3\textwidth,clip=true]{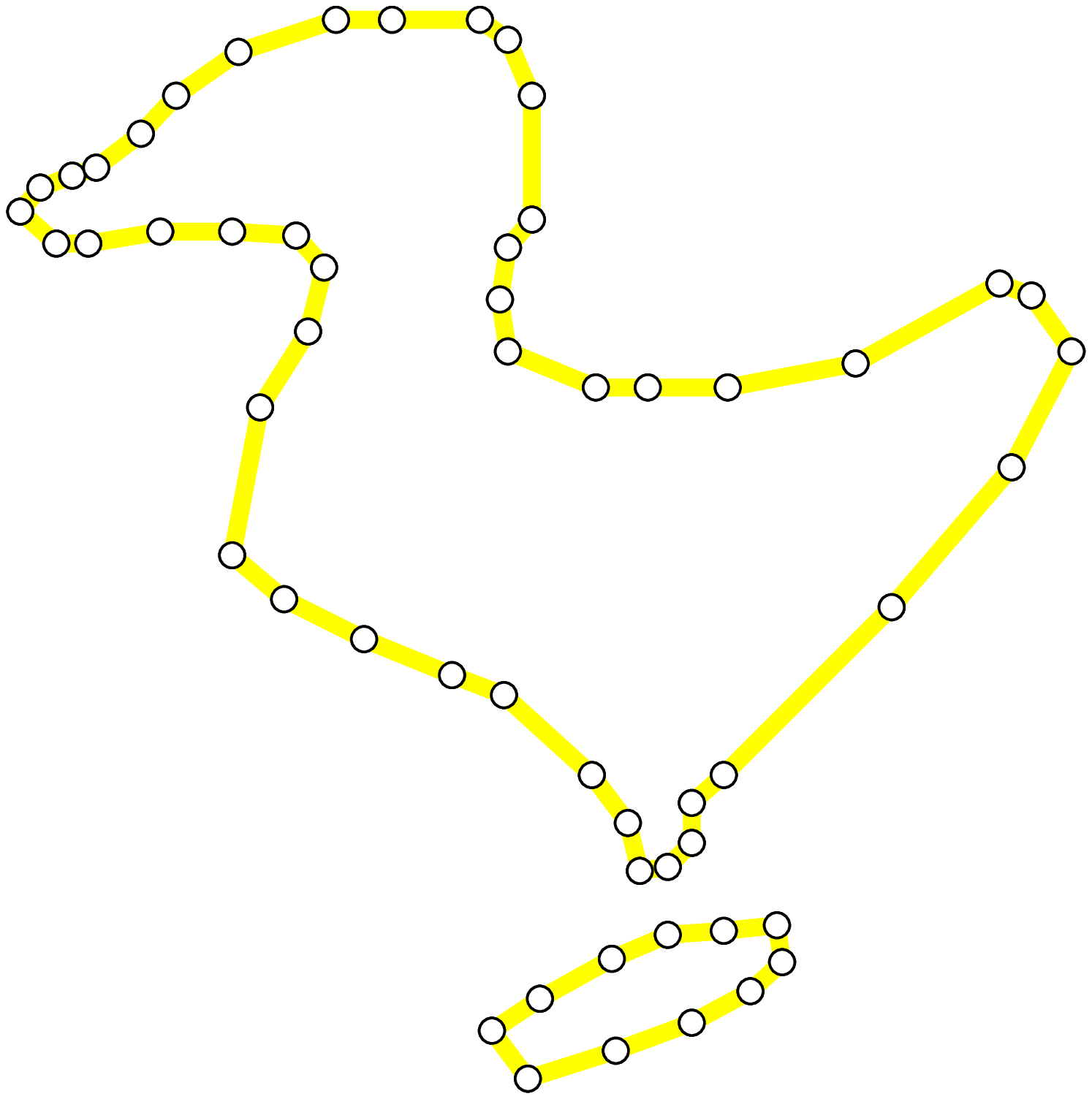}
  \caption{Planar curve(s) reconstructed from given points via the crust method of Amenta, Bern, and
    Eppstein~\cite{crust}.\label{fig:crust}}
\end{figure}

If a sufficiently well distributed finite set $S$ of points on a sufficiently smooth planar
curve~$K$ is given, then it is possible to obtain a polygonal reconstruction of~$K$.  Amenta, Bern,
and Eppstein~\cite{crust} obtained a curve reconstruction procedure via an iterated Voronoi diagram
computation.  This beautiful algorithm is implemented in \texttt{polymake}; see
Figure~\ref{fig:crust}.

\section{Selected Research Projects Using \texttt{polymake}}\label{sec:survey}
\noindent
We survey some projects for which \texttt{polymake} proved to be useful.

\subsection{Extremal Combinatorial Properties of Polytopes}

What originally inspired the development of the \texttt{polymake} system was the investigation of
combinatorial properties of convex polytopes and, in particular, the construction of interesting
examples.  Over the years \texttt{polymake} experiments helped to settle a number of previously open
questions. We list a few of them.

\subsubsection{Cubical Polytopes and Spheres} A \emph{cubical} polytope has only
(combinatorial) cubes as its faces.  Such objects naturally arise in computational geometry
(hexahedral meshes) as well as differential geometry (via normal crossing immersions of
hypersurfaces into spheres).

Algorithm complexity questions make it desirable to find out what the cubical polytopes with the
most number of faces with a given number of vertices are.  The main result
of~\cite{NeighborlyCubical} is the construction of so-called \emph{neighborly cubical polytopes}
which are more complex in this respect than previously expected.  A second paper on the same subject
gives deeper insight and it also describes more general constructions~\cite{math.CO/0503213}.

Schwarz and Ziegler~\cite{schwartz04:_const} construct a cubical $4$-polytope with an odd number of
facets.  This solves a previously open question.

\subsubsection{$f$-Vectors of $4$-Polytopes}

A classical result, due to Steinitz, characterizes all those triplets $(f_0,f_1,f_2)$ of natural
numbers for which there is a $3$-dimensional polytope with exactly $f_0$ vertices, $f_1$ edges, and
$f_2$ facets.  The corresponding question for higher dimensional polytopes is wide open.  In a
series of papers Ziegler et al. recently made progress as far as the understanding of such
\emph{$f$-vectors} of $4$-polytopes is concerned.  This progress, once again, is due to the
construction of special classes of $4$-polytopes, most of which were obtained with the help of
\texttt{polymake}.  See the survey~\cite{ziegler05:_convex_polyt} for an overview.

\subsection{Representation Theory of Groups}

Any representation $\nu:G\to\text{GL}(\RR^n)$ of a finite group~$G$ yields a \emph{representation
  polytope} $P(\nu)$ as the convex hull of the image in $\RR^{n^2}$.  This notion is introduced and
studied in Guralnick and Perkinson~\cite{math.CO/0503015}.  It turns out that combinatorial
properties of $P(\nu)$ are related to properties of the action of $\nu(G)$ on $\RR^n$.  For
instance, it is shown~\cite[Corollary 3.7]{math.CO/0503015} that if $\nu$ is a transitive
permutation representation then the diameter of the vertex-edge graph of $P(G)$ is bounded by
two. This is a far generalization of previous results on the Birkhoff polytopes due to Padberg and
Rao~\cite{MR0353997}.

\subsection{Gr\"obner fans}

A key task in algorithmic commutative algebra is the computation of Gr\"obner bases for a given
ideal in a polynomial ring.  From the theoretical as well from the practical viewpoint it is
essential that the resulting (reduced) Gr\"obner basis depends on a monomial ordering.  For deeper
structural investigations it is often useful to be able to understand the set of all Gr\"obner
bases for a given ideal.  The Gr\"obner fan is a polyhedral fan which gives a geometric structure to
this set.  Bahloul and Takayama~\cite{math.AG/0412044} used \texttt{polymake} in their quest to
extend the technique of Gr\"obner fans to ideals in other rings, in particular, the ring of formal
power series and the homogenized ring of analytic differential operators.

\subsection{Secondary Polytopes}

Triangulations of polytopes are interesting for various reasons.  For instance, they are
instrumental in solving certain systems of algebraic equations, via the theory of toric varieties
and sparse resultants.  The set of all triangulations of a fixed polytope~$P$ itself is endowed with
the structure of a (typically quite large) polytope, the \emph{secondary polytope} of~$P$.

Pfeifle and Rambau report~\cite{MR2011753} on an implementation for the construction of secondary
polytopes based on \texttt{TOPCOM}~\cite{topcom} and \texttt{polymake}.

\subsection{Enumeration of Small Triangulations of Manifolds}

The minimal number of vertices required to represent a given manifold as a simplicial complex is a
quite intricate topological invariant.  One approach to obtain this information for a large number
of manifolds is by complete enumeration of small triangulations.  This has been pursued by
Lutz~\cite{math.CO/0506372} and K\"ohler and Lutz~\cite{math.GT/0506520}. \texttt{polymake} was used
to compute simplicial homology groups and related invariants.

\subsection{Computational Biology}

A standard problem in computational biology is to deduce the optimal alignment of two given DNA
sequences.  Algorithms to solve this problem are known for some time.  However, the biologically
correct parameters which define what ``optimal'' means are usually not available.  This is a serious
obstacle on the way to biologically meaningful results.  Typically the computations have to rely on
estimates.  Thus it had been suggested to keep the parameters ``indeterminate'' and to compute with
the resulting algebraic expressions which then represent probabilities in the underlying statistical
models.  The impact of various choices of the parameter choices can then be studied in a
post-processing step.

While this naturally requires more complicated algorithms, it turns out by work of Pachter and
Sturmfels that methods from polyhedral geometry and \texttt{polymake}, in particular, can be employed.
This approach is comprehensively covered in the forthcoming book~\cite{ASCB}.

\section{Software Design}\label{sec:design}
\noindent
Since its first version from 1997 \texttt{polymake} was ---at least partially--- re-written several
times.  In spite of the many changes on the way, the core ideas always remained the same.  The first
goal was to have a flexible interface structure such that it is possible to interface to as many
existing polytope processing software components (developed by other people) as possible.  The
second goal was scalability in the sense that the system should be useful both for programmers and
mere users, and also both for students and expert scientists.

Feeling that one language is not enough for this, this resulted in an object-oriented hybrid design
based on the two programming languages \texttt{C++} and \texttt{Perl}.  The borderline is roughly
defined as follows: The \texttt{Perl} side takes care of all the object management and their
interfaces, while the \texttt{C++} half harbors the mathematical algorithms.  Information exchange
between both worlds is subject to a client-server scheme, \texttt{Perl} being the language of the
server.

\subsection{Open Objects}

Convex polytopes are represented in the system as a class of objects which are defined by an
extendible list of properties.  In the current distributed version there are already more than one
hundred of these properties defined; they range from the vertices (\texttt{VERTICES}) and facets
(\texttt{FACETS}) of a polytope to the list of Steiner points on all the faces
(\texttt{STEINER\_{}POINTS}) and the information whether or not the polytope is \texttt{SIMPLICIAL} or
\texttt{CUBICAL}.

The client perspective (on the \texttt{C++} side) is very restricted: The client asks the server for
properties of some polytope object and leaves it entirely to the server to decide how these should
be obtained.  The \texttt{Perl}-written server has a list of rules which specify how to compute
properties from the already known ones.  For instance, there is a rule which explains how the facets
can be computed for a polytope which was initially specified as the convex hull of finitely many
points; that is, there is a convex-hull algorithm rule which computes \texttt{FACETS} from
\texttt{POINTS}.  Actually, as in this case it is the fact, there may be several competing rules
computing the same.  It is the task of the server to compile admissible sequences of rules (via a
Dijkstra type algorithm for determining shortest weighted paths) to fulfill the user's (or the
client's) requests from the information initially given.

It is fundamental to the design that the set of rules as well as the list of properties known to the
system can be expanded and modified.  Moreover, the object management is abstract, too; this way it
is possible to define entirely new classes of objects and rule bases for them.  For instance,
simplicial complexes are objects different from polytopes (which actually includes pointed unbounded
polyhedra), while tight spans are specializations of polytope objects since they can be described as
the bounded subcomplexes of certain unbounded polyhedra; see Section~\ref{subsubsec:tightspan}.

\subsection{Scripting}

One way of using \texttt{polymake} is to generate large sets of polytopes and to filter them for
individual members with specific properties.  Such tasks are easily accomplished by often small
\texttt{Perl} scripts which make use of \texttt{polymake}'s object model.

As an example, the code below iterates through all the facets of a given polytope and offers a
visualization of the vertex-edge graphs of all the facets; on the way combinatorially equivalent
facets are detected and only one representative of each class is shown.

\begin{verbatim}
application 'polytope';

die "usage: polymake --script show_facets FILE\n" unless @ARGV;

my $p=load($ARGV[0]);
my @list=();

FACETS:
for (my $i=0; $i<$p->N_FACETS; ++$i) {
   my $facet=new Apps::polytope::RationalPolytope("facet #$i");
   Modules::client("facet", $facet, $p, $i, "-relabel");
   foreach my $other_facet (@list) {
      next FACETS if (check_iso($facet, $other_facet));
   }
   push @list, $facet;
}

static_javaview;
$_->VISUAL_GRAPH for @list;
\end{verbatim}

This script \texttt{show\_{}facets} is part of the distribution.

\subsection{Other Software}

There are also software packages which use \texttt{polymake} as a sub-system or have interfaces to
\texttt{polymake}.

\subsubsection{\texttt{convex}}

Franz has a \texttt{Maple} package for convex polytopes which is used, e.g., for homology
computations in toric varieties~\cite{convex}.

\subsubsection{\texttt{zerOne}}

L\"ubbecke implemented an algorithm for the enumeration of the vertices of a polytope with
$01$-coordinates given by its H-representation~\cite{zerOne}.

\subsubsection{\texttt{OpenXM}} Noro, Ohara, and Takayama~\cite{openxm} are the authors of
\texttt{OpenXM}, an infrastructure for mathematical communication.  This is used to interface
among several computer algebra systems.

\subsubsection{Electronic Geometry Models}

\begin{figure}[htbp]\centering
  \includegraphics[width=.5\textwidth]{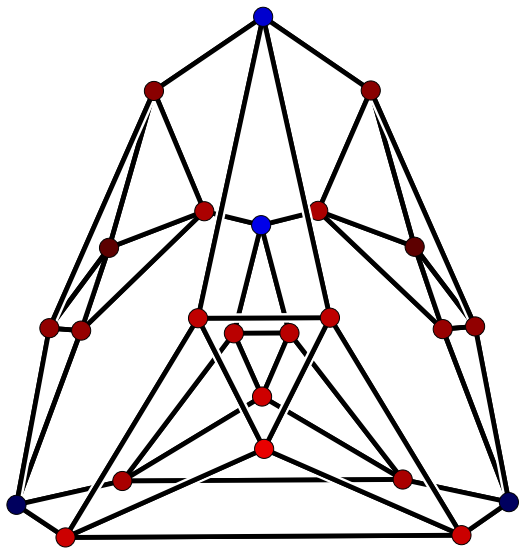}
  \caption{$3$-Polytope with four edges at each vertex, that can not be placed in convex position
    with all vertices on the surface of a sphere~\cite{2003.08.001}.\label{fig:eg-model}}
\end{figure}

Electronic Geometry Models is a refereed electronic journal for the publication of geometric
models~\cite{MMTCM:eg-models}.  Based on XML techniques these models are published in a standardized
way.  Some of these models describe polytopes, and these are given in \texttt{polymake} format; see
the example in Figure~\ref{fig:eg-model}.

The journal is freely accessible at \url{http://www.eg-models.de}.

\section{Technical Requirements}
\noindent
\texttt{polymake} can be used on UNIX systems only. It has been successfully tested on Linux, Sun
Solaris, FreeBSD, MacOS X, IBM AIX and Tru64 Unix.  Depending on the size of your objects
\texttt{polymake} can run on small machines with, say, 128 MB of RAM.  Only to compile the system
from the source code at least 1 GB of RAM is required.

Our website at \url{http://www.math.tu-berlin.de/polymake} offers the full source code as well as
several precompiled versions for download.

\bibliographystyle{amsplain}
\bibliography{main,polymake}

\end{document}